\documentclass[12pt]{article}
\usepackage{amssymb}
\usepackage{amsmath,amsfonts,amsthm}
\begin{document}
\newcommand{\inlim}{{lim}_{_{_{\!\!\!\!\!\!\!\!\!\!\!\longleftarrow}}}}
\newcommand{\drlim}{{lim}_{_{_{\!\!\!\!\!\!\!\!\!\!\!\longrightarrow}}}}
\newcommand{\fab}[2]{\langle {#1}_1,{#1}_2,\ldots ,{#1}_{#2} \rangle}
\newcommand{\abs}[1]{\:|#1|}    \newcommand{\Ker}{{\rm Ker\,}}
\newcommand{\ol}{\overline}     \newcommand{\Z}{ Z\!\!\!Z}
\newcommand{\C}{ l\!\!\!C}       \newcommand{\IH}{ I\!\!H}
\newcommand{\R}{ I\!\!R}         \newcommand{\N}{ I\!\!N}
\newcommand{\Q}{ l\!\!\!Q}       \newcommand{\notsubset}{ /\!\!\!\!\!\!\subset }
\newcommand{\id}{\mbox{id}}     \newcommand{\notdivide}{\!\not |\,}
\newcommand{\pf}{{\bf Proof: }}    \newcommand{\F}{ I\!\!F}
 \newcommand{\dg}{\mbox{deg}\,}
\newcommand{\st}{\; \vline\;}   \newcommand{\norm}[1]{\:\parallel #1 \parallel}
\newcommand{\orb}{\mbox{ \bf Orb}}  \newcommand{\stab}{\mbox{ \bf Stab}}
\newcommand{\M}{\mathcal{M}}      \newcommand{\rad}{\mbox{rad}}
    \newcommand{\ann}{\mbox{ann}}   \newcommand{\h}{\mbox{ht}}
    \newcommand{\Ass}{\mbox{Ass}}   \newcommand{\Supp}{\mbox{Supp}}
    \newcommand{\Spec}{\mbox{Spec}} \newcommand{\rk}{\mbox{rk }}
    \newcommand{\lr}{\mbox{little rank}}  \newcommand{\I}{\mbox{Im}}
    \newcommand{\Tor}{\mbox{Tor}} \newcommand{\n}{\mathcal{N}}
    \newcommand{\de}{ \mbox{depth}} \newcommand{\vcong}{\wr|}
     \newcommand{\td}{ \mbox{trdeg}}
\newcommand{\ba}{$$\begin{array}}\newcommand{\ea}{\end{array}$$}
\newcommand{\bea}{\begin{eqnarray*}}\newcommand{\eea}{\end{eqnarray*}}
\newcommand{\be}{\begin{equation}}\newcommand{\ee}{\end{equation}}
\newcommand{\bd}{\begin{definition}\bf} \newcommand{\ed}{\end{definition}}
\newcommand{\brs}{\begin{remarks}\rm}   \newcommand{\ers}{\end{remarks}}
\newcommand{\br}{\begin{remark}\rm}     \newcommand{\er}{\end{remark}}
\newcommand{\bt}{\begin{theorem}}       \newcommand{\et}{\end{theorem}}
\newcommand{\bl}{\begin{lemma}}         \newcommand{\el}{\end{lemma}}
\newcommand{\bco}{\begin{corollary}}    \newcommand{\eco}{\end{corollary}}
\newcommand{\bp}{\begin{proposition}}   \newcommand{\ep}{\end{proposition}}
\newcommand{\bo}{\begin{observation}\rm}\newcommand{\eo}{\end{observation}}
\newcommand{\bex}{\begin{examples}\rm}   \newcommand{\eex}{\end{examples}}
\newcommand{\bos}{\begin{observations}\rm}\newcommand{\eos}{\end{observations}}
\newcommand{\bx}{\begin{example}\rm}   \newcommand{\ex}{\end{example}}
\newcommand{\bexe}{\begin{exercise}\rm}   \newcommand{\eexe}{\end{exercise}}
\newcommand{\bpf}{\begin{proof}}
\newcommand{\epf}{\end{proof}}
      \newcommand{\K}{\mathcal{K}}

\newtheorem{definition}{Definition}[section]
\newtheorem{theorem}[definition]{Theorem}
\newtheorem{lemma}[definition]{Lemma}
\newtheorem{corollary}[definition]{Corollary}
\newtheorem{proposition}[definition]{Proposition}
\newtheorem{remarks}[definition]{Remarks}
\newtheorem{remark}[definition]{Remark}
\newtheorem{observation}[definition]{Observation}
\newtheorem{examples}[definition]{Examples}
\newtheorem{observations}[definition]{Observations}
\newtheorem{example}[definition]{Example}
\newtheorem{exercise}[definition]{Exercise}

\title{ON POWER STABLE IDEALS}
\author{Pramod K. Sharma\\ e--mail:
 pksharma1944@yahoo.com\\
 School Of Mathematics, Vigyan Bhawan, Khandwa Road,\\ INDORE--452
017, INDIA.}
\date{}
\maketitle
                  \begin{center}
                  \section*{ABSTRACT}
                  \end{center}
     We  define the notion of a power stable ideal in a polynomial ring $ R[X]$ over an
     integral domain $ R $. It is proved that a maximal ideal $\chi$ $  M $ in $ R[X]$ is
     power stable if and only if $ P^t $ is $ P$- primary for all $ t\geq 1 $
      for the prime ideal $ P = M \cap R $. Using this we prove that for a
       Hilbert domain $R$ any radical ideal in $R[X]$ which is a finite intersection
      G-ideals is power stable. Further, we prove that if $ R $ is  a Noetherian
      integral domain of dimension 1 then any radical ideal in $ R[X] $ is power stable.
      Finally, it is proved that if  every ideal in $ R[X]$ is power stable then $ R $ is
       a field.
        \section{INTRODUCTION }
              All rings are commutative with identity $(\neq 0)$. For a subset $S$ of a
              ring $R$, $id(S)$ shall denote the ideal of $R$ generated by $S$ and for
              an ideal $J$ in $R,\widehat R_J$ shall denote the $J-$adic completion of the
              ring $R$. If $ I $ is an ideal in $R[X]$, then for any $a(X)\in R[X]$,
               $\bar{a}(X)$ will denote the image of $a(X)$ in $R[X]/I$.
          In [2], we define an ideal $I$ in a ring $R$ almost prime ideal if for all $a,b\in R,
         ab\in I-I^2$ either $a\in I$ or $b\in I$. While trying to prove that all ideals
          in $\Z[X]$ are almost prime we required that for an ideal $I$ in $\Z[X]$,
       $I\cap \Z= n{\Z} $ implies $I^2\cap \Z = n^2\Z$. This, however, was not true.
        The property seems interesting in itself and is the basis of our definition
        of power stable ideal in a polynomial ring $R[X]$. For an integral domain $R$,
        an ideal $I$ in $R[X]$ is called power stable if $I^t\cap R = (I\cap R)^t$ for
        all $t\geq 1$. Hereafter, a ring $R$ shall
        always denote an integral domain. In this note we initiate the study of power
        stable ideals. We prove that a maximal ideal $ M $ in $ R[X] $ is power stable
         if and only if for the prime ideal $P = M\cap R $,
       $P^t$ is $ P$ - primary for all $t\geqslant 1$ (Theorem 3.7). This
       result is used to prove that if $R$ is a Hilbert domain then any radical ideal
       in $R[X]$ which is a finite intersection of  G-ideals is power stable (Theorem 3.10).
        Further, it is proved that if $R$ is a Noetherian domain of  dimension 1 then any
         radical ideal in $R[X]$ is power stable (Theorem 3.11). We also prove that if
         every ideal in $R[X]$ is power stable then $R$ is a field (Theorem 3.14).

        \section{Observations On Definition}
          We define:
        \bd\label{1.1} Let $R$ be an integral domain. An ideal $I$ in $R[X]$ is called
        power stable ideal if for all $t\geq 1$, $I^{t}\cap R = (I\cap R)^t$.
         \ed
          \bx \label{1.2}Any principal ideal in $R[X]$ is power stable.\ex
              \bx\label{1.3}For any ideal $I$ of $R$ the ideal $I[X]$ is power

              stable.\ex
                \bx \label{1.4}If $I$ is a power stable ideal of $R[X]$ then $I^{t}$ is power
                 stable for all $ t\geq 1 $. \ex

       First we make some general observations.
       \bl Let $I$ be an ideal in $R[X]$, and $I\cap R=J$. Then $I$
       is power stable if and only if the natural homomorphism $$\phi :
      Gr_JR\longrightarrow Gr_IR[X]$$ is monomorphism of graded rings. \el
      \bpf If $I$ is power stable, then $I^n\cap R=J^n$ for all $n\geq 0$.
      Hence \ba{rcl} J^n\cap I^{n+1} & = & I^n\cap R \cap I^{n+1} \\ & = &
       I^{n+1} \cap R = J^{n+1} \ea Therefore $\phi$ is a monomorphism.
       Conversely, let $\phi$ be a monomorphism. Then \be \label{111.1} J^n\cap
       I^{n+1}=J^{n+1} \mbox{ for all } n\geq 0 \ee We shall prove that
       $I^n\cap R= J^n$ for all $n\geq 1$ by induction on $n$. Since $\phi$ is
       monomorphism the statement is clear for $n=1$. Let $n\geq 2$. By induction
       assumption, \ba{rl} & I^{n-1}\cap R=J^{n-1} \\ \Rightarrow &
        I^n\cap (I^{n-1}\cap R) =  I^{n}\cap J^{n-1}\\ \Rightarrow & I^n\cap R = J^n  \\ \ea  by equation
       \ref{111.1}. Hence the result follows. \epf
        \bl Let $I$ be an ideal in
       $R[X]$ and $I\cap R=J$. Then $I$ is power stable if and only if the
       natural map $$\widehat{R}_J \longrightarrow \widehat{R[X]}_I$$
       is a monomorphism. \el
       \bpf First of all, note that $\widehat{R}_J = \inlim R/J^n$ and
       $\widehat{R[X]}_I = \inlim R[X]/I^n$. If $I$ is power stable then
        $I^n\cap R = J^n$ for all $n\geq 1$. Hence the natural map
        $R/J^n\stackrel {\alpha_n}{\longrightarrow}  R[X]/I^n$ is a monomorphism,
         and the diagram :\ba{rcl} R/J^{n+1}& \stackrel {\alpha_{n+1}}{\longrightarrow}&
         R[X]/I^{n+1}\\ \downarrow & & \downarrow\\ R/J^n & \stackrel
         {\alpha_n}{\longrightarrow}&
         R[X]/I^n \ea is commutative where vertical maps are quotient maps. This set up
          induces a natural monomorphism $ \widehat{R}_J = \inlim R/J^n \longrightarrow
         \widehat{R[X]}_I = \inlim R[X]/I^n$ .  On the other hand, if $I$ is not
         necessarily power stable, $\alpha_n$ still exists and the diagram above is
         commutative. Thus the set up induces a homomorphism from $\widehat{R}_J$ to
          $\widehat{R[X]}_I$. Now it is easy to see that if this homomorphism is a
           monomorphism then $I$ is power stable.                  \epf
       \section{Main Results}
        \bl Let $I$ and $J$ be power stable ideals in $R[X]$ such that $I\cap R$
       and $J\cap R$ are co-maximal. Then $I\cap J$ is power stable.\el
       \bpf  Let $t\geq 1$. Then \ba{rl} (I\cap J)^t\cap R & \subset  (I^t\cap R)\cap ( J^t\cap
       R)  \\
       & = (I\cap R)^t\cap (J\cap R)^t \mbox { since } I,J \mbox{ are power stable. }  \\
          & = ((I\cap R)\cap (J\cap R))^t \mbox { since } I\cap R \mbox{ and } J\cap R
           \mbox { are co-maximal. }
          \\&  = (I\cap J\cap R)^t .\ea
           Further, it is clear that $ (I\cap J\cap R)^t \subset (I\cap J)^t \cap R $.
            Therefore  $(I\cap J\cap R)^t = (I\cap J)^t \cap R $. Consequently $I\cap J $
              is power stable. \epf
        \bl An ideal $I$ in $R[X]$ is power stable if and only if $I_{\it P}$
        is power stable in $R_{\it P}[X]$ for all $\it P\in \Spec\ R .$
        \el
         \bpf  If $I$ is power stable then for any prime ideal $P$ in $R$, $I_P$ is power
          stable since localization commutes with intersection and powers. Further, we
           always have $(I\cap R)^n \subset I^n\cap R$. If $I_P$ is power stable
           for a prime ideal $P$ in $R$ then ${(I_P)}^n \cap R_P = (I_P\cap R_P)^n $.
            Therefore $((I^n\cap R)/(I\cap R)^n)_P = 0$. Hence if $I_P$ is power stable for all
           primes in $R$ then clearly  $I^n\cap R = (I\cap R)^n$ i.e., $I$ is
           power stable. \epf

              \bt Let $R$ be a principal ideal domain. Let for an ideal $I$ in
              $R[X]$, $I\cap R = Rd$. If the image of $I$
              in $R/(d)[X]$ is generated by a regular element then $I$ is power
              stable ideal.
              \et
              \bpf  By assumption on $I$, $I= id(d,h(X)) $ where image of $h(X)$ in
              $R/(d)[X]$ is a regular element. We shall prove by induction on $t$
              that $I^t\cap R = (I\cap R)^t$  for all $t\geq 1.$ Let $t\geq 2,$
              and let $I^s\cap R = (I\cap R)^s$  for all
              $s \leq t-1$. If $I^t \cap R = eR$, then
               \ba {rrll}
                 & eR & = & I^{t} \cap R \subseteq I^{t-1} \cap R = d^{t-1}
              R\\
               \Rightarrow & e & = & d^{t-1} k, \mbox { where }  k\in R
               \ea
               Now, as $e\in I^t$, we have
               \ba {rrll}
                 & d^{t-1}k  & = & \displaystyle \sum_{i+j=t, i\geq 1} d^i h^j(X) a_{ij}(X) + h^t
                (X)a(X) \mbox{ for some }  a_{ij}(X), a(X) \in R[X].\\
                \Rightarrow &  \overline{h}(X)^t \overline{a}(X) & = & 0  \mbox { in }
                R/(d)[X]\\
                \Rightarrow & \overline{a}(X) & = & 0  \mbox{ in } R/(d)[X], \mbox{ since }
                \overline{h}(X)
                \mbox{ is regular in  } R/(d)[X]\\
                \Rightarrow & a(X) & = & d b(X) \hspace{.1in } ( b(X)\in R[X] )\\
                \Rightarrow & d^{t-1} k & = & \sum d^i h^j(X) a_{ij}(X)+ d h^t(X) b(X)\\
                \Rightarrow & d^{t-2}k & = &  \sum d^{i-1} h^j(X) a_{ij}(X)+  h^t(X) b(X)\\
                 \Rightarrow & d^{t-2}k & \in & I^{t-1}\cap R=d^{t-1}R \\ \Rightarrow &
                 k & \in & dR \\
                 \Rightarrow & e & = & d^tk_1\in d^tR \hspace{.1in}\mbox{ for }k = dk_1.
                 \ea
                   Hence $ I^t\cap R=d^tR $, and the result follows.\epf
                 \bco Let $R$ be a principal ideal domain. If $P$ is a prime ideal in $R[X]$, then $P$ is power stable.
                 \eco
                 \bpf We have either $P\cap R = (0)$ or $Rp$ where $p$ is a prime element in
                 $R$. As any non zero prime ideal in $R$ is  maximal ideal, the proof
                 is immediate from the theorem.\epf

       \bt\label{111.4} Let for an ideal $I$ in $R[X],$ $I=id(J,\,f(X))$ where $J$
       is an ideal in $R$ and $f(X)\in R[X]$ is a monic polynomial of degree
       greater than or equal to 1. Then $I$ is a power stable. \et \bpf To
       prove the result we have to show that $I^t\cap R=J^t$ for all $t\geq
       1$. Let us first consider the case $t=1$. If $\lambda \in I\cap R$,
       then we can write $$\lambda = a(X)+f(X)h(X)$$ where $a(X)\in J[X]$.
       Reading off this equation in $R/J[X]$, we get
       \ba{rl} & \bar{\lambda}=\bar{f}(X)\bar{h}(X) \\ \Rightarrow &
       \bar{\lambda}=0, \mbox{since $\bar{f}(X)$ is  monic polynomial
       of degree $\geq 1$ in $R/J[X]$ . }\\ \Rightarrow & \lambda \in J. \ea
        Hence $I\cap R=J$, and the result follows for $t=1$. \\ Now, let $t>1$ and
        $\lambda \in I^t\cap R$. Then we can write \be\label{111.2}
       \lambda = a(X)+f(X)h(X)+f^t(X)c(X) \ee where $a(X)\in J^t[X]$ and
       $h(x)\in J[X]$. As in the case $t=1$, reading off this equation in
       $R/J[X]$, we conclude, $c(X)\in J[X]$ and $\lambda\in J$. Hence we
       can write \be \label{111.3} \lambda = a(X)+f(X)b(X) \ee where
       $a(X)\in J^t[X]$ and $b(X)\in J[X]$ such that no coefficient of
       $b(X)$ is in $J^t$. Now, if $a_0$ is leading coefficient of $a(X)$ and $b_0$ is
       leading coefficient of $b(X)$, then $a_0+b_0=0$. This implies
       $b_0\in J^t$. A contradiction to our assumption. Hence $\lambda\in J^t$.
        Therefore $I^t\cap R=J^t=(I\cap R)^t$.This completes the proof.\epf
         \bco\label{111.5} Suppose for an ideal $I$ in $ R[X],$
          $I\cap R= M$ is a maximal ideal. Then $I$ is power stable. \eco
         \bpf If $I= M[X]$, then clearly $I$ is power stable. Hence, let $I\neq M[X]$.
         Now, as $ I\cap R= M, $ it is clear that $I=id(M,f(X))$
         where $f(X)\in R[X]$ is a monic polynomial of
       degree greater than or equal to 1. Therefore the result follows from the
       theorem . \epf \bt A maximal ideal $M$ in $R[X]$ is power stable if and
       only if for $P=M\cap R$, $P^{(t)}=P^t$ for all $t\geq 1$ i.e., $P^t$
       is $P$-primary for all $t\geq 1$. \et \bpf Let $M$ be power stable.
       As $M^t$ is $M$-primary for all $t\geq 1$, $M^t\cap R=P^t$ is $M\cap
       R=P$-primary for all $t\geq 1$. Conversely, let $P^t$ be
       $P$-primary for all $t\geq 1$. If $P=(0)$, there is nothing to
       prove. Hence, let $P\neq(0)$. Then as \ba{rl} & M_P\cap R_P[X]=PR_P
       \\ \Rightarrow & M_P \mbox{ is power stable by Corollary
       \ref{111.5}. } \\ \Rightarrow & M^t_P\cap R_P=P^tR_P \\ \Rightarrow
       & M^t\cap R\subseteq M^t_P\cap R\subseteq P^tR_P\cap R= P^{(t)} \\
       \Rightarrow & M^t\cap R\subseteq P^{(t)}=P^t \\ \Rightarrow & M^t\cap
       R=P^t=(M\cap R)^t \ea  Thus the result is proved. \epf \br (1) In the
       reverse part of above result it is not used that $M$ is maximal.
       Thus if $P$ is a prime ideal in $R[X]$ and for $p=P\cap R$,
       $p^t=p^{(t)}$ for all $t\geq 1$, then $P$ is power stable. Further,
       note that if $P$ is a power stable prime ideal in $R[X]$, then for
       $p=P\cap R$, $p^t$ need not be $p$-primary for all $t\geq 1$. This
       is clear since for any $p\in Spec\ R$, $P=p[X]$ is power stable
       prime in $R[X]$. Thus if $p^t$ is not $p$-primary, we get the
       required example. \\ (2) If $R$ is a Hilbert domain then any
       maximal ideal in $R[X]$ is power stable. In particular, if $K$ is a
       field then any maximal ideal in  the polynomial ring $K[X_1, \cdots X_n]$
        is power stable.  \er
        We now give an example to show that, in general, a maximal
        ideal in $R[X]$ need not be power stable.
         In view of Theorem 3.7, it suffices to give a $G-$ ideal $ P $
         in $ R $ for which $ P^n $ is not $ P- $ primary for some $n\geq 1.$
       The example below was suggested by Melvin Hochster.
        \bx  Let $K$ be a field and $ Y,Z,W $ be algebraically independent elements
         over $ K .$ For an algebraically independent element $ T $ over $ K ,$
          consider the $ K-$ algebra homomorphism :
                     $$ \phi :K[[Y,Z,W]]\longrightarrow K[[T]] $$
         such that $ \phi(W)= T^3$, $\phi(Y)= T^4$ and $\phi(Z)=T^5$. Then kernel of
          $\phi$ is the prime ideal $P$= id$(f,g,h)$ where $ f=(W^3-YZ),g=(Y^2-WZ) $
          and $ h=(Z^2-W^2Y) $. It is easy to see that $P$ is a G-ideal. Further, we
          have $ f^2 - gh = W Q(W,Y,Z)\in P^2.$  Clearly $ W\not\in P$ and it is easy
          to check that $Q \not \in P^2$. Thus $ P^2 $ is not $ P- $primary.
           Hence for the integral domain $R = K[[Y,Z,W]]$, there is a
             maximal ideal $M$ in the polynomial ring $R[X]$ such that $M\cap R= P$
              and $M$ is not power stable.\ex

        \bt\label{111.6} Let $R$ be a Hilbert ring. Then any radical ideal $I$ in
        $R[X]$ which is a finite intersection of of G-ideals, is power stable. \et
       \bpf By [3, Theorem 31], $R[X]$ is a Hilbert ring. Hence all $G$-ideals in
        $R[X]$ are maximal. Now, by our assumption on $I$, $$I=M_1\cap \cdots \cap
       M_k,$$ where $M_i$'s are distinct G- ideals in $R[X]$. As $R$
       is Hilbert ring, $M_i\cap R=m_i$ is maximal ideal in $R$. Now note
       that for any $t\geq 1$, \ba{rcl} I^t & = & M_1^t\cap \cdots \cap
       M_k^t\\ \Rightarrow I^t\cap R & = & M_1^t\cap \cdots \cap
       M_k^t \cap R \\ & = & m_1^t\cap m_2^t\cap \cdots \cap m_k^t \ea
        since, by Corollary 3.6, every maximal ideal in $R[X]$ is power stable. $$\Rightarrow
       I^t\cap R=m_{i_1}^tm_{i_2}^t\cdots m_{i_l}^t$$ where
       $m_{i_1},\,m_{i_2}\cdots m_{i_l}$ are all distinct maximal ideals
       in the set $\{ m_i\,|\,\,1\leq i\leq k \}$. Thus clearly $I$ is
       power stable. \epf \bt Let $R$ be a Noetherian domain of dimension 1. Then
       any radical ideal in $R[X]$ is power stable. \et \bpf  If $I\cap R=(0)$, there is
       nothing to prove. Hence, assume $I\cap R = J(\neq 0)$. Since $I$ is
        radical ideal in $R[X]$, $J$ is a radical ideal in $R$. As $R$ is
        Noetherian domain of dimension 1, we have $$ J = M_1\cap \cdots \cap M_k$$  where $M_i$'s
        are maximal ideals in $R$. Thus it is clear that for any prime ideal $P$ in
         $R$ either $J_P = R_P$ or $J_P = PR_P$. Therefore, since $ I_P\cap R_P = J_P$
         for every prime ideal $P$ in $R$,  by Corollary 3.6,
         $I_P$ is power stable for any prime ideal $P$ in $R$. Hence by Lemma 3.2,
          $I$ is power stable. \epf

         We shall now show that in case $R$ is of dimension 1, a non-radical ideal
        in $R[X]$ need not be power stable. In fact, we shall give an  example of a
       primary ideal in $R[X]$ which is not power stable where $R$ is a principal
       ideal domain. This example was given by Melvin Hochster( personal
       communication) for $ R = \Z$. We learnt this via Stephen McAdam.
                   \bx\label{p123} Let $R$ be a P.I.D. and $p$ be a prime in
                     $R$. Then $I =id(X^2-p, X^3)$ is not a power stable
                     ideal in $R[X]$.
                                    \ex
                         \bpf  Step 1. $I\cap R = Rp^2$.\\
                                 We have\\
                                 \ba{lrrl}
                                        & X^3- X(X^2-p) & = &  pX\in I \\
                                           \Rightarrow & pX^2 & \in &I \\
                                           \Rightarrow & pX^2 -p(X^2-p)& = & p^2\in
                                           I
                                           \ea
              If $I\cap R=Rd$, then $d$ divides $p^2$. We, now note :\\ $(i)$ $d\neq
              1$ \\ If $ d=1 $, then \ba{ll} & 1 = (X^2-p)a(X)+X^3b(X)
              \hspace{.2in} (a(X), b(X) \in R[X]) \\ \Rightarrow & 1 = -pa(0) \\
              \Rightarrow & p \mbox{ is a unit } \ea This is absurd. Thus $d\neq 1$.\\
              $(ii)$ $d\neq p$ \\ If $d=p$, then we can write $$p=(X^2-p)a(X)+X^3b(X)
              \qquad (a(X), b(X) \in R[X]) $$ Putting $X=\sqrt p$ in the above equation,
               we get $$ p= (p\sqrt p)  b(\sqrt p)$$ Clearly
              $$b(\sqrt p)=c+d\sqrt p, \mbox{ for some } c, d \in R.$$ Thus \ba{lrrl}
              & p & = & p\sqrt p (c+d\sqrt p) \\ \Rightarrow & 1 & = & c\sqrt p + dp
              \\ \Rightarrow & c^2p & = & 1+d^2p^2-2dp \\ \Rightarrow & 1 & = &
              p(c^2-d^2p+2d) \ea This implies $p$ is a unit, which is not true.
               Thus step 1 is proved. \\ Step 2. $p^3 \in I^2\cap R$ \\ We have $$I^2=id(X^4-2pX^2+p^2,
              X^6, X^5-pX^3)$$ Now \ba{lrrl} & X^6-X(X^5-pX^3) & = & pX^4 \in I^2 \\
              \Rightarrow & pX^4-p(X^4-2pX^2+p^2) & = & 2p^2X^2-p^3 \in I^2 \ea As
              seen in step 1, $pX\in I$. Therefore $p^2X^2 \in I^2$, and hence $$
              p^3=2p^2X^2-(2p^2X^2-p^3) \in I^2$$ This proves step 2. \\ By step 1 and
              step 2 it is immediate that $I^2\cap R\neq (I\cap R)^2$. \epf \br In the
              above example radical of $I$ is $id(X, p)$, a maximal ideal in $R[X]$.
              Thus $I$ is a primary ideal. Hence primary ideals need not be power
              stable. \er
               In the end we note:
              \bt Let $R$ be a Noetherian integral domain. If every ideal in $R[X]$ is
              power stable, then $R$ is a field. \et \bpf Let $ P $ be a prime ideal of
              height $1$ in $R$. By our assumption, every ideal in $R_P[X]$ is power
              stable. Thus, to prove our result, we can assume $R$ is local ring of
              dimension $1$. Let $(R, M)$ be a local ring of dimension $1$. By
              assumption, for any $\lambda\in M-M^2$, $J=id(X^2-\lambda, \lambda X)$ is
              power stable. Put $I=J\cap R$. If $a\in J\cap R=I$, then \ba{rl} &
              a=f(X)(X^2-\lambda)+g(X)\lambda X \,\,\,(f(X),\,g(X)\in R[X]) \\ \Rightarrow &
              a=-f(0)\lambda \,\,\mbox{ putting } X=0. \\ \Rightarrow & I\subset
              R\lambda \\ \Rightarrow & I=I_1\lambda \,\,\mbox{ where } \,\,
              I_1=\{b \in R\,|\,\, b\lambda \in I\}.\ea We, now, consider two
              cases.\\ Case $1$: $I_1\subset M$. \\ Let us note that $ J^2 = id( X^4
              - 2\lambda X^2 + \lambda^2 , \lambda^2 X^2, \lambda X^3 - \lambda^2 X) .$
              Therefore, as $\lambda^3 = \lambda( X^4 - 2\lambda X^2 +\lambda^2) +
              \lambda^2 X^2 -  X(\lambda X^3 - \lambda ^2 X)),$ we have
                \ba{rl}
               & \lambda^3\in J^2\cap R=I_1^2\lambda^2 \\ \Rightarrow &
               \lambda^3=b\lambda^2\,(b\in I_1^2) \\ \Rightarrow & \lambda = b \in
               I_1^2\subset M^2 \ea This is not true by choice of $\lambda$. Therefore
               case $1$ is not possible. \\ Case $2$: $I_1=R$. \\In this case $J^2\cap
               R=R\lambda^2$ i.e., $\lambda^2\in J^2$. Hence \ba{rl} & \lambda^2 \in
               id(X^4-2\lambda X^2+\lambda^2,\,\lambda^2 X^2,\,\lambda X^3-\lambda^2 X^) \\ \Rightarrow &
               \lambda^2= (X^4-2\lambda X^2+\lambda^2)a(X)+\lambda^2 X^2 b(X) +
                (\lambda X^3  -\lambda^2 X)c(X) \\ \Rightarrow & \lambda^2\equiv
               X^4 a(X)(\mbox{mod}\,\lambda ) \ea i.e., $\lambda^2$ is
               a multiple of $X^4$ in $R/(\lambda)[X]$. This is absurd. Hence $R$ is
               a field and the result follows. \epf
              It would be interesting to know the answer to the following :\\
                \bf{ Question.} Let $R$ be a Noetherian integral domain of
                dimension 1. Does there exist a characterization of power stable
                ideals in $R[X]$?
              \begin{center}\section*{ACKNOWLEDGEMENT}\end{center}
              I am thankful to Stephen McAdam for some useful e-mail exchanges. I am
              also thankful to Melvin Hochster for the examples and more so for
              his very prompt responses to all my questions. Finally, I express my
              thanks to the referee for his suggestions on presentation and for
              expressing suspicion about the proof  of the original version of
              Theorem 3.10 .

              \begin{center}
  {\bf REFERENCES}
  \end{center}
\begin{enumerate}
\item  M.F.Atiyah, I.G. Macdonald ; Introduction to Commutative Algebra, Addison- Wesley Publ. Co., 1969.

\item  S.M. Bhatwadekar, Pramod K. Sharma; Unique factorization and birth of Almost
Primes,Communications in Algebra (To appear)

\item  Irving Kaplansky, Commutative rings, The university of Chicago Press, 1974.
\end{enumerate}
\end{document}